\documentclass[12pt,a4paper]{article}
\usepackage[T1]{fontenc}
\usepackage{apalike}
\usepackage{enumerate}
\usepackage{setspace}
\usepackage{graphicx}
\usepackage{subfloat}
\usepackage{endfloat}
\usepackage{fancyhdr}
\usepackage{amsmath}
\usepackage{amssymb}
\usepackage{amsfonts}
\usepackage{color}

    \newcommand{\InRefFig}[1]{Figure~\ref{#1}}
    
    \newcommand{\RefFig}[1]{Fig.~\ref{#1}}
    
    \newcommand{\RefEq}[1]{Eq.~(\ref{#1})}
    \newcommand{\RefEqs}[1]{Eqs.~(\ref{#1})}
    \newcommand{\RefsEqs}[2]{Eqs.~(\ref{#1}) and (\ref{#2})}
    \newcommand{\Ref}[1]{(\ref{#1})}
    \newcommand{\InRefSec}[1]{Section~\ref{#1}}
    
    \newcommand{\RefSec}[1]{Sec.~\ref{#1}}
    \newcommand{\RefSecs}[1]{Secs.~\ref{#1}}

\makeatletter
\renewcommand{\listoffigures}{\section*{Figure Captions}\@starttoc{lof}}
\makeatother

\begin{document}

\title{Harmonic analysis of oscillators through standard numerical continuation tools}

\markboth{IJBC, Bizzarri, Linaro, Oldeman, Storace: Harmonic
analysis of oscillators by numerical continuation}{}

%
%
\singlespacing
\begin{center}
\verb+F. Bizzarri, D. Linaro, B. Oldeman, M. Storace+\\
\verb+Harmonic analysis of oscillators through standard numerical continuation tools+\\
\verb+International Journal of Bifurcation and Chaos+\\
\verb+to appear on vol. 20, n. 12, Copyright @ 2010 World Scientific Publishing Co.+\\
\verb+With permission from World Scientific Publishing Co. Pte. Ltd.+
\vspace{2cm}

{\LARGE \bf Harmonic analysis of oscillators through standard
numerical continuation tools \vspace{1.5cm}}

{\large Federico Bizzarri,$^a$ Daniele Linaro,$^a$ Bart Oldeman,$^b$
Marco Storace$^{a,*}$
\vspace{0.75cm}
\\ $^a$ Biophysical and Electronic Engineering Department
\\ University of Genoa, Via Opera Pia 11a, Genova, Italy
\vspace{0.25cm}
\\ $^b$ Department of Computer Science and Software Engineering
\\ Concordia University, 1455 De Maisonneuve Blvd. West,
\\ Montreal, Quebec, Canada
\vspace{0.25cm}
\\ E-mail: {\texttt{marco.storace@unige.it}}
\vspace{1.5cm}}
\end{center}

%
%
\begin{center}{\bf Abstract}\end{center}
\noindent In this paper, we describe a numerical continuation method
that enables harmonic analysis of nonlinear periodic oscillators.
This method is formulated as a boundary value problem that can be
readily implemented by resorting to a standard continuation package
- without modification - such as AUTO, which we used. Our technique
works for any kind of oscillator, including electronic, mechanical
and biochemical systems. We provide two case studies. The first
study concerns itself with the autonomous electronic oscillator
known as the Colpitts oscillator, and the second one with a
nonlinear damped oscillator, a non-autonomous mechanical oscillator.
As shown in the case studies, the proposed technique can aid both
the analysis and the design of the oscillators, by following curves
for which a certain constraint, related to harmonic analysis, is
fulfilled.

\singlespacing

\section{Introduction}
Continuation methods provide an efficient tool for analyzing systems
of nonlinear algebraic equations whose solutions form a
one-dimensional continuum. When dealing with periodic solutions of
systems of ordinary differential equations (ODEs), we continue
solutions by solving a boundary value problem (BVP). We can either
solve this BVP explicitly, when possible, or implicitly, by using
numerical continuation tools.

The BVPs that are used to continue ``standard'' objects, such as
equilibrium points, periodic, homoclinic, and heteroclinic orbits and
their bifurcations in dynamical systems are well-known (see, for
instance, \cite{Kuznetsov:2004} for an overview). However, it is often
advantageous to formulate new BVPs to continue ``non-standard''
objects, such as invariant manifolds \cite{Doedel:2006}, slow
manifolds \cite{Desroches:2008}, and coherent structures such as
spiral waves and other defects in oscillatory media
\cite{Bordyugov:2007,Champneys:2007}.

In a recent paper \cite{Cochelin:2009}, a combination of harmonic
analysis and continuation techniques (based on the harmonic balance
method) was proposed. The main result is a new numerical continuation
tool that provides standard results of continuation analysis after a
preliminary reformulation of the problem in terms of harmonic balance.

In this paper, we show how to define a novel BVP that enables harmonic
analysis by using standard numerical continuation tools.  The proposed
BVP allows, in the spirit of
\cite{Doedel:2006,Desroches:2008,Bordyugov:2007,Champneys:2007}, for
``non-standard'' continuations focused on selected harmonic components
of a solution without an \textit{ad hoc} simulator like the one
described in \cite{Cochelin:2009}.

The first application in this paper deals with the \textit{design} of
an autonomous electronic circuit (the Colpitts oscillator). Its aim is
to obtain parameter charts that can be immediately understood by
designers of electronic oscillators. In particular, these charts show
how the limit cycle corresponding to a periodic regime changes in a
familiar circuit parameter plane. The guidelines for the designer are
provided, for instance, by curves corresponding to solutions with
fixed amplitudes of some harmonic components or with a fixed ratio of
two harmonic components.

The second application entails the harmonic \textit{analysis} of a
non-autonomous mechanical oscillator (nonlinear damped oscillator).
In particular, we use the proposed BVP to analyze the so-called jump
phenomenon \cite{Schmidt:1986}.

The main advantages of the proposed approach are the following:
\begin{itemize}
  \item the BVP is formulated in such a way that
  the public-domain software package AUTO-07P \cite{Auto07PMan}
  can solve it;
  \item compared to simulations, computation times are generally lower,
  since numerical continuation packages operate directly on system
  invariant sets;
  \item the proposed procedure is reasonably easy to use
  by those who are familiar with the analysis of nonlinear
  dynamical systems.
\end{itemize}
On the other hand, one needs to take care of the following aspects:
\begin{itemize}
  \item for complex systems, such as a realistic radio-frequency
  electronic oscillator, the procedure is effective only
  if preceded by a modeling phase where one defines a suitable
  model of the oscillator \cite{Bizzarri:2009} ---
  such a model should be as simple as possible, but able to capture the
  essential features of the system;

  \item those who are not familiar with numerical continuation
    packages require a preliminary training.

\end{itemize}

This paper is organized as follows. \InRefSec{sec:basic} summarizes
the basic elements that the proposed technique is based on. The BVP
formulation is described in \RefSec{sec:setup}, whereas
\RefSecs{sec:ex1} and \ref{sec:ex2} are devoted to two case studies.
Some conclusions are drawn in \RefSec{sec:concl}.

\section{Basic elements}\label{sec:basic}
Let the following system of ODEs describe an oscillator, which can be
either autonomous or non-autonomous with periodic forcing:
\begin{equation}
\dot{x} = \frac{dx}{dt} = g(x,t;p), \quad x,g \in \mathbb{R}^n, \; p
\in \mathbb{R}^q, \; t \in \mathbb{R}.
\label{eq:ODE}
\end{equation}
\noindent
When this oscillator reaches a periodic regime, it produces signals of
generic period $T$, that is, frequency $f = \frac{1}{T}$ and angular
frequency $\omega = 2\pi f$. The Fourier series expansion of each
state variable is
\begin{equation}
x_j(t) = a_{0j}+\sum_{k=1}^\infty \left[ a_{kj} \sin(k \omega t) +
b_{kj} \cos(k \omega t) \right],
\label{eq:fourier}
\end{equation}
\noindent where the index $j \in \{1,\ldots,n\}$ selects the state
variable, $a_{0j}$ is the mean value of $x_j(t)$ over $T$ and we
obtain the other coefficients by projecting $x_j(t)$ on the
corresponding basis functions
\begin{equation}
\begin{split}
a_{kj} &= \frac{2}{T} \int_t^{t+T} x_j(\tau) \sin(k \omega \tau) d\tau\\
b_{kj} &= \frac{2}{T} \int_t^{t+T} x_j(\tau) \cos(k \omega \tau) d\tau\\
\end{split}
\label{eq:coeffourier}
\end{equation}

In general, continuation methods show how system invariants (for
example, equilibrium points or limit cycles) depend on one or more
control parameters. One of the key elements of continuation theory is
that invariant sets and their bifurcations are revealed when so-called
\textit{test functions} equal zero. Many test functions are included
in the most diffused numerical continuation packages
\cite{Auto07PMan,Dhooge:2003}, but, of course, user-defined test
functions can be added. For instance, in this paper we define test
functions of the kind $f(S_a,S_b,T,K_\text{REF})$, where $S_a$ and
$S_b$ denote a subset of $n_a$ and $n_b$ Fourier coefficients
$\{a_{kj}\}$ and $\{b_{kj}\}$, respectively, and $K_\text{REF}$ is a
constant reference value.

We can use this formulation to solve different kinds of problems.
Firstly we can continue the amplitude of a harmonic component with
respect to a single parameter. However, of further interest are
continuations with respect to two parameters, adding a further
constraint. For instance, iso-harmonic, iso-ratio, and iso-energy
continuations are possible. For iso-harmonic continuations we fix the
amplitude of a harmonic component, for iso-ratio continuations the
ratio between amplitudes of different harmonics, and for iso-energy
continuations a sum of squares of a limited number of harmonic
amplitudes representing almost the whole power spectrum of the
analyzed signal.

In the next section, we set up the BVP that is solved for
these continuations.

\section{Setup}\label{sec:setup}
We start by deriving a simplified model of the oscillator given by a
small system of ODEs, as is usual when dealing with both the analysis
and synthesis of a dynamical system.  We call this system the
\textit{original system}.  By assuming some reasonable modeling
hypotheses we obtain these equations, which are expressed in terms of
state variables.  It is advantageous, but not compulsory, to normalize
these equations and shift the origin of the normalized state space to
a ``significant'' equilibrium point. Such an equilibrium is stable for
some parameter configuration, and, by varying one parameter, undergoes
a supercritical Andronov-Hopf bifurcation, which marks the appearance
of a family of asymptotically stable periodic solutions evolving
around the (unstable) equilibrium.

\subsection{The original system}
For autonomous oscillators, the original system is simply the ODE
system modeling the oscillator, given by $\dot x = g(x(t);p)$.  For a
non-autonomous oscillator with periodic forcing, we can obtain an
equivalent autonomous oscillator by adding a nonlinear oscillator with
the desired periodic forcing as one of the solution components (see,
for instance, the AUTO-07P demo \texttt{frc}
\cite{Auto07PMan,Alexander:1990}). In particular, for a sinusoidal
forcing we can use the following secondary oscillator (very close to
the normal form of the supercritical Andronov-Hopf bifurcation):

\begin{equation} \label{eq:auxosc}
\begin{split}
\dot v &= \alpha v + \beta w + v (v^2+w^2)\\
\dot w &= -\beta v + \alpha w + w (v^2+w^2),\\
\end{split}
\end{equation}
which for $\alpha<0$ asymptotically converges to the origin and for
$\alpha>0$ has the asymptotically stable solution $v = \sin(\beta
t)$, $w = \cos(\beta t)$. For instance, if the first state variable
of the original system obeys the differential equation $\dot{x}_1 =
g_1(x(t);p) + c \cos(\omega t)$ and the state of the original system
is $x = [x_1, \ldots, x_n]$, then in the corresponding autonomous
system the equation for the first state variable is given by
$\dot{x}_1 = g_1(x(t);p) + c w$, where $\beta = \omega$ and the
state vector is redefined as $x = [x_1, \ldots, x_n,v,w]$. Hence, if
the original system is a non-autonomous oscillator with periodic
forcing, we can also recast it to the autonomous system $\dot{x} =
g(x(t);p)$.

We can analyze equilibria, their Hopf bifurcations and emanating
limit cycles in the possibly recast original system. The original
system, together with some information from the limit cycle close to
the Hopf bifurcation, is then used to construct and initialize the
full continuation system defining the BVP.

\subsection{Initialization}
As an initial solution for the BVP problem, we could take a
pre-computed periodic orbit, carry out a signal analysis to find the
coefficients of its Fourier expansion, and substitute that into the
system. Another way (which we describe here) is to use standard
continuation techniques provided by AUTO, following an equilibrium
point that undergoes an Andronov-Hopf bifurcation.

In the non-autonomous case with sinusoidal forcing, we can easily find
this bifurcation by varying the parameter $\alpha$ in
\RefEqs{eq:auxosc} across zero.

The limit cycle that emanates from the Andronov-Hopf bifurcation can
then be continued a small distance away from the bifurcation, where it
will still have the approximate form
\begin{equation} \label{eq:simplefourier}
x(t)=A \sin \omega t + B \cos \omega t.
\end{equation}
Here $\omega$ is the purely imaginary part of the corresponding
eigenvalue of the equilibrium, which we can obtain using the period
$T$ that AUTO provides: $\omega = 2\pi / T$. Now the Fourier
coefficients can be trivially derived, comparing
\RefsEqs{eq:fourier}{eq:simplefourier}:
\begin{equation} \label{eq:fouriercoeff}
\begin{aligned}
a_{0j} &= 0, \\
a_{1j} &= A=\dot{x_j}(0)/\omega, \quad \quad & b_{1j} &= B=x_j(0),\\
a_{kj} &= 0, \quad \quad & b_{kj}&=0, \quad k>1, \quad j=1,\ldots,n\\
\end{aligned}
\end{equation}

\subsection{The full continuation system} The full continuation
system is given by the following equations:
\begin{itemize}
\item
Non-autonomous differential equations:
\begin{equation}
\begin{split}
\dot x &= T g(x(t),t; p), \quad x, g \in \mathbb{R}^n, \quad p \in \mathbb{R}^q, \quad t, T \in \mathbb{R}\\
\dot t &= 1\\
\end{split}
\label{eq:full}
\end{equation}
Here the possibly recast original system \RefEq{eq:ODE} is rescaled so
that a $T$-periodic solution of \RefEq{eq:ODE} is a 1-periodic
solution for the first equation in \RefEqs{eq:full}. The non-periodic
equation $\dot t=1$ is added to make the system look autonomous to
continuation software.

\item
Boundary conditions that define a periodic orbit on the $t$-interval
$[0,1]$:
\begin{equation}
\begin{split}
  x_j(0) &= x_j(1), \quad j = 1,\ldots,n\\
  t(0) &= 0\\
\end{split}
\end{equation}
\item
Integral conditions:
\begin{equation}
\begin{split}
  \int_0^1 x(t) \dot{x}_{\text{old}}(t) dt &= 0\\
  \int_0^1 (2 x_j(t) \sin(2\pi k t) - a_{kj}) dt &= 0 \quad \text{for any} \; a_{kj} \in S_a\\
  \int_0^1 (2 x_j(t) \cos(2\pi k t) - b_{kj}) dt &= 0 \quad \text{for any} \; b_{kj} \in S_b\\
  \int_0^1 f(S_a,S_b,T,K_\text{REF}) dt &= 0\\
\end{split}
\end{equation}
Here the first condition is the standard integral phase condition
\cite{Kuznetsov:2004}, where $x_{\text{old}}(t)$ denotes the previous
point on a continuation branch, the second and third conditions
compute or fix the Fourier coefficients $\{a_{kj}\}\in S_a$ and
$\{b_{kj}\}\in S_b$, and the fourth condition computes or fixes
$K_\text{REF}$ given the Fourier coefficients in $S_a$ and $S_b$.
\end{itemize}

This gives us a system of $n+1$ ordinary differential equations, $n+1$
boundary conditions and $n_a+n_b+2$ integral conditions.  Adding a
standard pseudo-arclength condition, this needs to be offset by
$n_a+n_b+3$ continuation parameters.

A basic choice for the parameters is to continue a periodic orbit in
one of the system parameters, the period $T$, and the $n_a+n_b+1$
values in $S_a \cup S_b \cup \{K_\text{REF}\}$. Note that in this case
these last $n_a+n_b+1$ values are effectively measured through
integral conditions, where explicit test functions such as
\begin{equation}
a_{kj} =  \int_0^1 (2 x_j(t) \sin(2\pi k t) ) dt
\end{equation}
would suffice.

However, the integral conditions become more powerful if we like to
keep something fixed: for instance, by fixing $K_\text{REF}$ and
freeing up one more system parameter we can continue iso-amplitude
curves. Moreover, in existing continuation software it is easier and
arguably more elegant, if a little more computationally expensive, to
stick to one full system with the same integral conditions, than to
switch between test functions and equivalent integral conditions.

\section{Case study 1: an autonomous electronic oscillator} \label{sec:ex1}
Increasing demand, in modern communication systems, of
high-performance low-power radio-frequency circuits is driving the
development of accurate simulation tools at the design stage. The
simulation is particularly critical in the case of analog circuits
(free-running or voltage-controlled oscillators), because of time
consumption, and in the case of mixed-signal analog-digital circuits
(frequency synthesizers or phase-locked loops), due to the modelling
difficulties. In all cases, the analysis is usually non-trivial and
performance verification requires extensive simulation. This is even
more important when the main goal is to find how changes in circuit
parameters (for example, the amplitude or frequency of an input
generator, or a linear-element value) affect circuit performance.

The main research lines in this field are twofold.  Firstly there is
the development of algorithms that speed up and make circuit
simulations more reliable \cite{Brambilla:2005,Brambilla:2008}.
Secondly there is the study of methods that, starting from a
simplified version of the designed circuit, aim to provide design
criteria avoiding brute-force simulations of complex integrated
circuits. Most of these methods are based on well-established theories
such as harmonic balance, Volterra series (for weakly nonlinear
oscillators), and bifurcation analysis of nonlinear dynamical
systems. Harmonic balance is used for an accurate determination of
nonlinear-circuit operation bands \cite{Rizzoli:1988,Suarez:1998} and,
associated with other methods, for calculating the periodic response
of nonlinear dynamical circuits \cite{Buonomo:2003}. Volterra series
are used for the analysis of nearly sinusoidal nonlinear oscillators
\cite{Hu:1989,Huang:1994}. Bifurcation analysis is used to optimize
the design of oscillators \cite{Maggioetal:1999} and the locking range
of injection-locked frequency dividers \cite{Ghahramani:2007}. A
combination of harmonic-balance simulators and bifurcation control is
exploited to obtain bifurcation control in microwave circuits, thus
presetting the operation bands of complex circuits, such as
synchronized and voltage-controlled oscillators and frequency dividers
\cite{Collado:2005}.  This combination is also used to obtain robust
and efficient oscillator analysis techniques (see, for example,
\cite{Bonani:1999,Genesio:1993,Gourary:2000} and references therein).

\begin{figure}[h!!]
\centering{\includegraphics{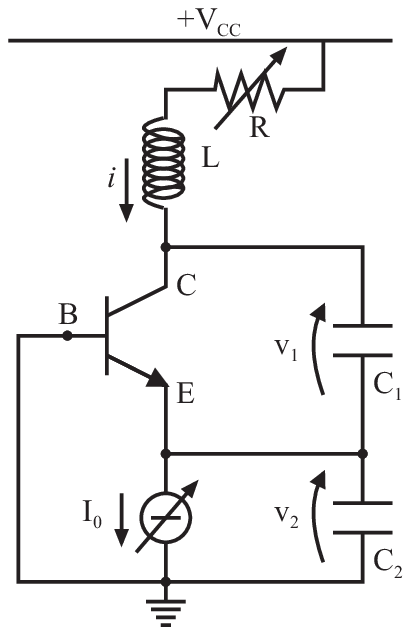}}

\caption{The considered Colpitts oscillator. $V_{CC}=3$ V,
$C_1=C_2=1$ $\mu$F, $L=10$ mH.}
\label{fig:colpitts}
\end{figure}
Our case study is the Colpitts oscillator, shown in
\RefFig{fig:colpitts}, whose dynamics were analyzed in detail in
\cite{Maggioetal:1999,DeFeoetal:2000,DeFeoMaggio:2003}. The
following simplified model, that we can easily substitute into the
system defined in Section~\ref{sec:setup}, adequately describes this
oscillator:
\begin{equation} \label{eq:Colpitts}
\left\{
\begin{aligned}
\dot{x} &= \frac{G}{Q(1-\gamma)}\bigl[-\alpha_F \left(e^{-y}-1
\right)+z\bigr]
\\
\dot{y} &= \frac{G}{Q
\gamma}\bigl[(1-\alpha_F)\left(e^{-y}-1\right)+z\bigr]
\\
\dot{z} &= -\frac{Q \gamma(1-\gamma)}{G}\bigl(x +
y\bigr)-\frac{1}{Q}z
\end{aligned} \right. ~,
\end{equation}
where the system parameters are initialized as follows:
\begin{equation} \label{eq:Colpittsnormpar}
\begin{gathered}
Q = \frac{\omega_0 L}{R} = 0.8; \quad \gamma=\frac{C_2}{C_1+C_2};
\quad G=\frac{I_0 L}{V_T R (C_1
+ C_2)} = 2;\\
C_1 = C_2 = 10^{-6} [F]; \quad L = 10^{-3} [H].
\end{gathered}
\end{equation}

Moreover, we assume $\alpha_F=1$, that is, we assume the CB
short-circuit forward current gain of the BJT transistor to be
ideal. The variable $z$ denotes the inductor current normalized with
respect to $I_0$.  The variables $x$ and $y$ are the voltages across
$C_1$ and $C_2$ normalized with respect to $V_T = 25.9$ mV (that is,
the thermal voltage at room temperature). Time is normalized with
respect to $T_0 = \sqrt{L C_1 C_2 / (C_1 + C_2)}$.

Continuing the equilibrium at $0$ as the parameter $G$ goes downwards
from $2$ we find a Hopf bifurcation at $G=1$. For $G$ slightly greater
than $1$ there then exists a limit cycle for which the Fourier
coefficients $a_{12}$, $b_{12}$, and $K_\text{REF}$ are given by (see
\RefEqs{eq:fouriercoeff})
\begin{equation}
\begin{split}
  a_{12} &= \dot y(0)/(2\pi/T) = \frac{G z(0)}{Q \gamma}/(2\pi/T),\\
  b_{12} &= y(0),\\
  K_\text{REF} &= \sqrt{a_{12}^2+b_{12}^2}.\\
\end{split}
\end{equation}

We focus on the state variable $y$ since it corresponds to the
voltage used as output of the oscillator. Moreover, we focus on its
first harmonic component since we want the generated oscillation to
be nearly sinusoidal. Thus we need to verify that higher order
harmonics have small amplitudes when compared with the first
harmonic component.

We then follow the periodic orbits as a boundary value problem in $(G,
T, a_{12},$ $b_{12}, K_\text{REF})$ with user-defined labels for the
following values of $K_\text{REF}$: $1$, $4$, $7$, $10$, $13$, $16$,
$19$, $22$, $25$. Using these values as starting points we can find
iso-harmonic curves in $(G,\gamma)$-parameter plane for fixed values
of $K_\text{REF}$ by continuing in $(Q, G, T, a_{12}, b_{12})$.

\begin{figure*}[h!!!]
\centering{\includegraphics[scale=1.0]{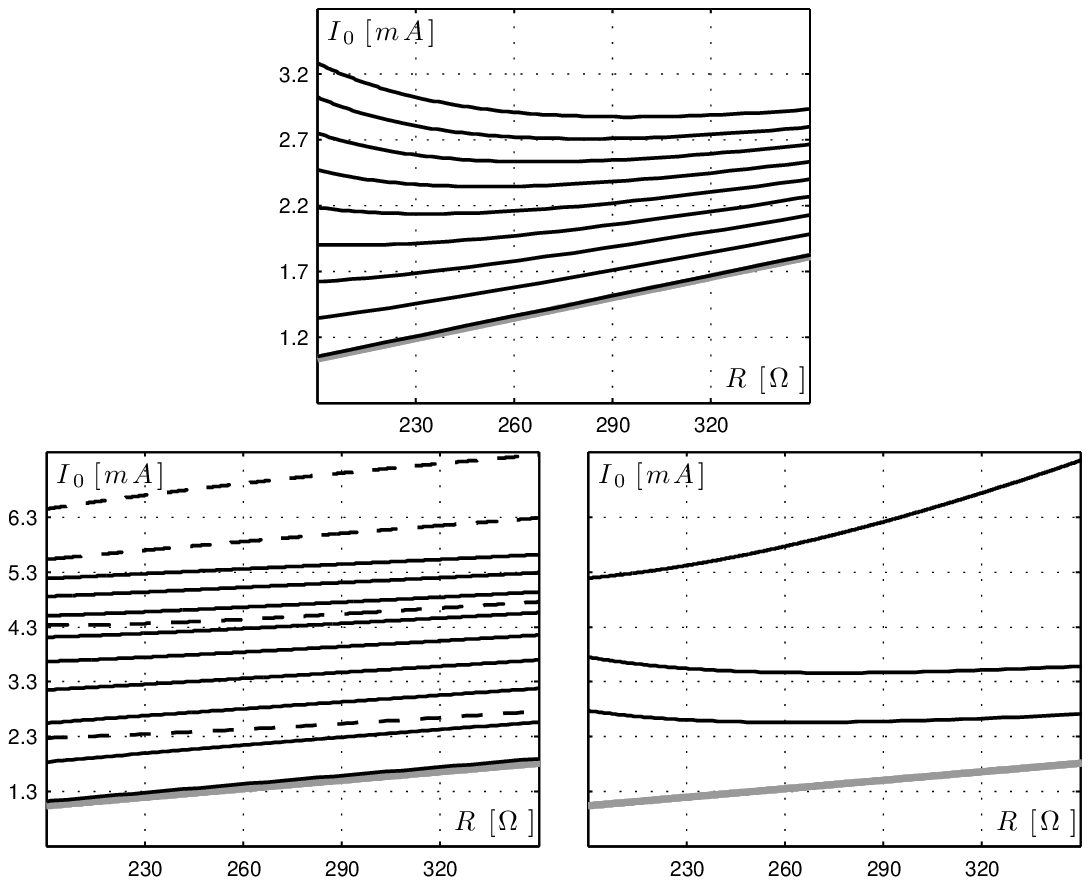}}

\caption{Charts for the Colpitts oscillator on the parameter plane
$(R,I_0)$. The grey line is the Hopf bifurcation curve in all
panels. Iso-period curves (upper panel), iso-$A_1$ (solid) and
iso-$A_2$ (dashed) curves (lower-left panel), and iso-ratio
$\frac{A_1}{A_2}$ (lower-right panel).}
\label{fig:colpittsres1}
\end{figure*}

We carried out the continuation in the parameter space $(Q, G, T,
a_{12}, b_{12})$, but show the results in \RefFig{fig:colpittsres1} on
the plane $(R,I_0)$, according to the definitions of $Q$ and $G$, to
simplify the interpretation of the results from a circuit point of
view. We excluded the parameter space region characterized by the
presence of complex dynamics \cite{Maggioetal:1999,DeFeoetal:2000}
from our analysis and focused on the region characterized by the
presence of a stable ``simple'' harmonic cycle.

In the upper panel, along the black dashed curves the period
(normalized with respect to $T_0$) is constant: $6.3$ for the lower
curve, and $8.7$ for the upper curve, with a step of $0.3$.  In
the lower-left panel, along the black solid curves the amplitude $A_1$
(normalized with respect to $V_T$) of the first harmonic component of
$y$ is constant ($1$ for the lower curve, $25$ for the upper curve,
step of $3$), whereas along the black dashed curves the amplitude
$A_2$ (normalized with respect to $V_T$) of the second harmonic
component of $y$ is constant ($1$, $5$, $10$, and $15$ from the lower
curve to the upper one).

In the left-right panel, we follow the periodic orbits as a boundary
value problem in $(G, T, a_{12}, b_{12}, a_{22}, b_{22},
K_\text{REF})$. Along the black solid curves the ratio
$\frac{A_1}{A_2} \; (= K_\text{REF})$ is constant: $3$ for the upper
curve, $5$ for the lower curve, with a step of $1$. The
interpretation of these results is straightforward: for instance, if
we want to fix both circuit parameters to have an oscillation
frequency $f=T_0/6.3$, and ensure low values of $R$ to increase the
quality factor $Q$, we can just properly adjust $I_0$ along the lower
curve in the upper panel of \RefFig{fig:colpittsres1}. The results
shown in the lower panels provide useful information if a periodic
signal whose first harmonic component is predominant with respect to
the second one is of interest.

Of course, we can obtain other iso-curves besides those shown in
\RefFig{fig:colpittsres1}. For example, we may chose values of
$I_0$ and $R$ ensuring a desired ratio $\frac{A_p}{A_q}$ for given
$p$ and $q$.

\begin{figure}[t!!]
\centering{\includegraphics[width=0.8\columnwidth]{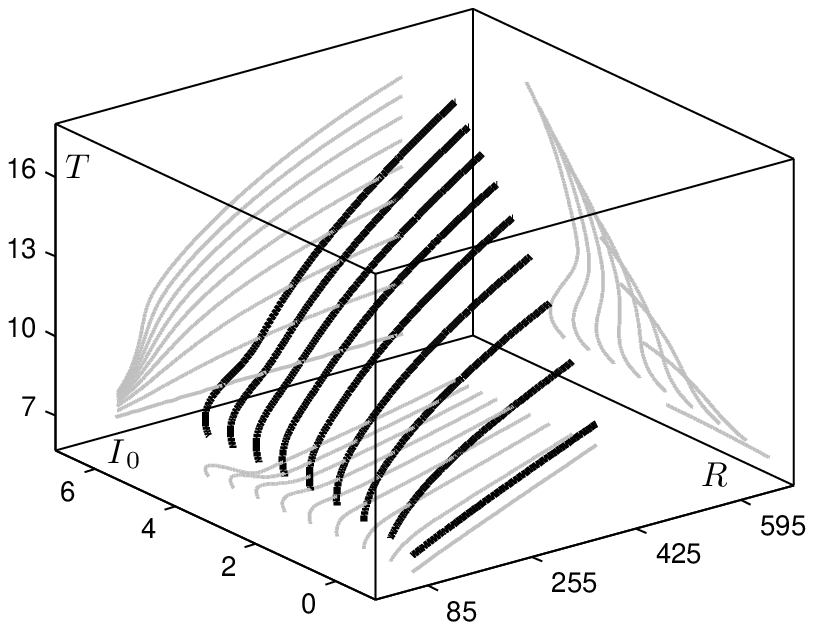}}

\caption{Iso-harmonic curves in the parameter subspace $(R,I_0,T)$.}
\label{fig:colp3d}
\end{figure}

\InRefFig{fig:colp3d} shows the obtained iso-harmonic (with constant
$A_1$) curves in the parameter subspace $(R,I_0,T)$. Of course, the
projection of the curves on the plane $(R,I_0)$ gives the black
solid curves displayed in the lower-left panel of
\RefFig{fig:colpittsres1}.

\section{Case study 2: a non-autonomous mechanical oscillator} \label{sec:ex2}
The nonlinear damped oscillator is a model widely used to represent
shock absorbers \cite{Lang:2006} and is given by the following
system (with dimensionless variables and parameters):
\begin{equation} \label{eq:ndo}
\begin{split}
 \dot{x} &= \frac{A}{m}\cos(\omega t) - \frac{1}{m} (c_1 x + c_2 x^2 + c_3 x^3 + k y) \\
 \dot{y} &= x \\
\end{split}
\end{equation}
where $m=240$, $c_1 = 296$, $c_2 = 3000$, $c_3 = 800$, and $k = 240 (4
\pi)^2$. The reference angular frequency of the system is $\omega_0
= \sqrt{k/m}=4\pi$.

We can easily recast the non-autonomous system as an autonomous
system by adding the auxiliary oscillator \Ref{eq:auxosc} to
\RefEqs{eq:ndo}. For $\alpha < 0$, the resulting system has a stable
equilibrium at the origin, which undergoes a Hopf bifurcation for
$\alpha = 0$. For small positive $\alpha$ there then exists a limit
cycle for which the Fourier coefficients $a_{12}$, $b_{12}$, and
$K_\text{REF}$ are approximately given by, and determined
numerically as (see \RefEqs{eq:fouriercoeff})
\begin{equation}
\begin{split}
  a_{12} &= \dot y(0)/(2\pi/T) = x(0)/(2\pi/T),\\
  b_{12} &= y(0),\\
  K_\text{REF} &= \sqrt{a_{12}^2+b_{12}^2}.\\
\end{split}
\end{equation}

We focus on the state variable $y$ since it represents the position
of the mechanical oscillator.

Once the initialization is completed, we work with the full
continuation system and continue the limit cycle with respect to
$(\alpha, T, a_{11}, b_{11}, K_\text{REF})$ until we reach $\alpha =
1$, corresponding to the correct sinusoidal forcing.

\begin{figure*}[t!!]
\centering{\includegraphics[width=1.0\columnwidth]{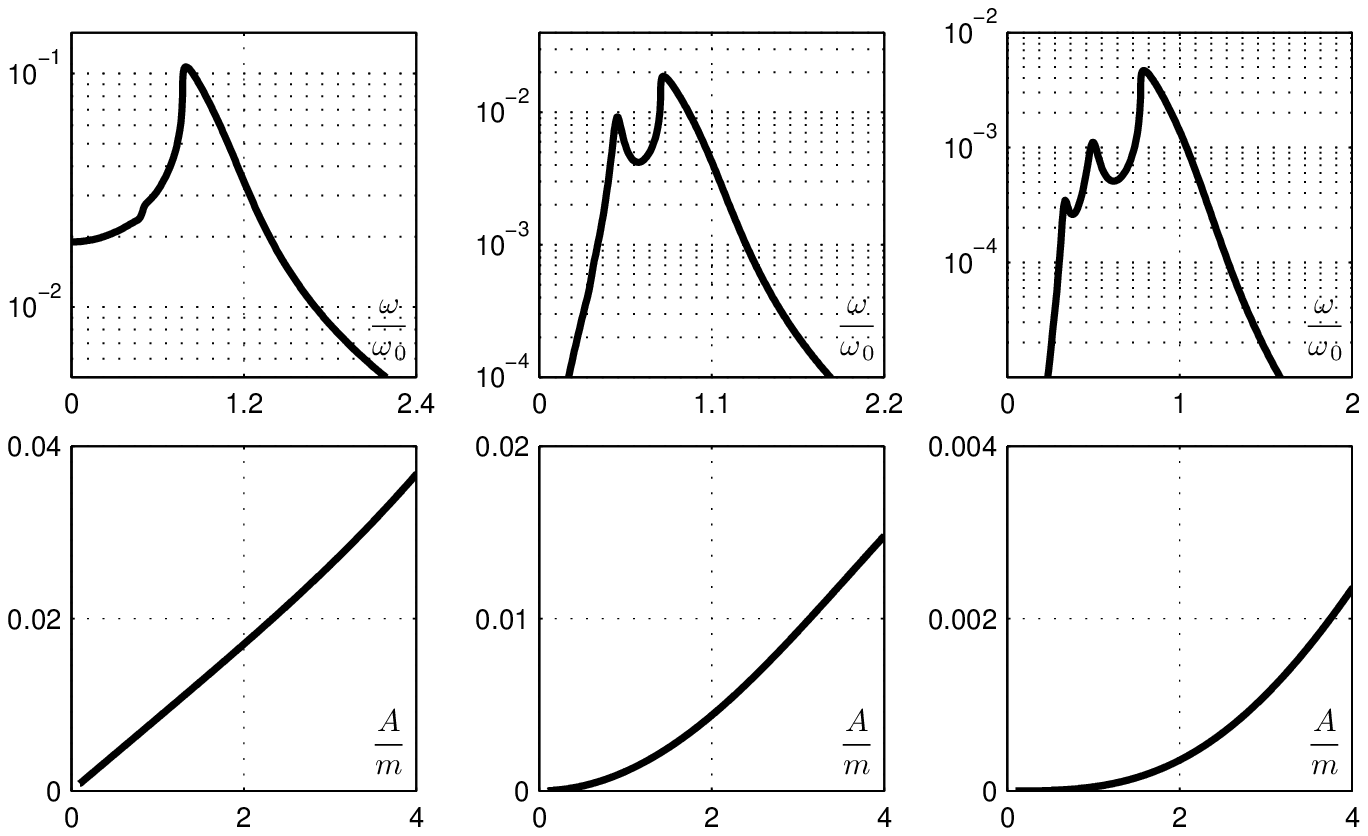}}
\caption{Amplitudes of the first (first column), second (second
column) and third (third column) harmonics with respect to
$\omega/\omega_0$ (upper panels) and to $A/m$ (lower panels).}
\label{fig:ndoresults}
\end{figure*}

In this case, we monitor the amplitude of the first harmonic with
respect to one of the parameters $\omega$ and $A$ (properly
normalized), besides $(T,a_{11},b_{11},K_\text{REF})$, to provide
evidence for the presence of a ``jump phenomenon''. The ``jump
phenomenon'' is a characteristic feature of many nonlinear
oscillators, where the response amplitude changes suddenly at some
critical value of the excitation frequency \cite{Schmidt:1986}. Many
widely used methods such as the Harmonic Balance (HBM) and Nonlinear
Output Frequency Response Function (NOFRF) \cite{Peng:2008} barely
capture this phenomenon. \InRefFig{fig:ndoresults} shows the
amplitudes of the first three harmonics as functions of the
normalized frequency $\omega/\omega_0$ (upper panels) and of the
normalized amplitude $A/m$ (lower panels). The results perfectly
match the reference diagrams reported in \cite{Peng:2008}.

We remark that the proposed continuation method does not require any
approximation, thus providing excellent accuracy. Moreover, it is
based on tools (such as AUTO) that are reliable and widely tested.

Finally, by performing continuations similar to those shown for the
Colpitts oscillator, we can also obtain simulations pertaining to the
design of mechanical oscillators, that is, we obtain the values of
parameters that ensure a desired behavior.

\section{Conclusions} \label{sec:concl}
We proposed a technique, based on numerical continuation, that enables
harmonic analysis of a given nonlinear oscillator without resorting to
any approximation. Moreover, a designer can choose some oscillator
parameters to obtain a desired behavior, by analysing some of the
curves that are obtained by this technique.  More realistically, since
the model is an approximation of the real system, the proposed method
provides at least reference values of the bifurcation parameters. More
accurate simulations focused on restricted portions of the parameter
space can then refine these values.

The main advantage of this technique is that it enables the analysis
of even relatively complex oscillators (both autonomous and forced) by
using software tools that are reliable and optimized. This makes the
development of \textit{ad hoc} software for this kind of analysis
unnecessary and ensures an excellent accuracy of the results.

The main limit of this technique is that it requires a thorough
knowledge of continuation methods and/or software packages for
numerical continuation. Moreover, for oscillators forced by
non-sinusoidal periodic signals, it can be non-trivial to define the
auxiliary equations needed to make the system autonomous. In these
cases, it may be more convenient to use the results of a signal
analysis carried out on a pre-computed periodic orbit. The
coefficients of the thus obtained Fourier expansion can be substituted
into the system, which can then be used as an initial solution for the
BVP problem.

\section*{Acknowledgments}
The authors acknowledge helpful discussions with Eusebius Doedel.


\begin{thebibliography}{99}
%
\bibitem[Alexander {\em et~al.\/}, 1990]{Alexander:1990}
    {{J.C.} Alexander, {E.J.} Doedel, \& {H.G.} Othmer},
    ``On the resonance structure in a forced excitable system,''
    {\it SIAM J. Appl. Math.}, vol. 50, pp. 1373--1418, 1990.
%
\bibitem[Bizzarri {\em et~al.\/}, 2009]{Bizzarri:2009}
    {F.~Bizzarri, A.~Brambilla, D.~Linaro, \& M.~Storace},
    ``Continuation analysis of a phase/quadrature electronic oscillator,''
    {\it Journal of Circuits, Systems and Computers},
    special issue on {\em Advances in oscillator analysis and design}, in press.
%
\bibitem[Bonani \& Gilli, 1999]{Bonani:1999}
    {F.~Bonani \& M.~Gilli},
    ``Analysis of stability and bifurcations of limit cycles in {C}hua's circuit through the harmonic-balance approach,''
    {\it {IEEE} Transactions on Circuits and Systems {I}: Fundamental Theory and Applications}, vol. 46, pp. 881--890, 1999.
%
\bibitem[Bordyugov \& Engel, 2007]{Bordyugov:2007}
    {G.~Bordyugov \& H.~Engel},
    ``Continuation of spiral waves,''
    {\it Physica D}, vol. 228, pp. 49--58, 2007.
%
\bibitem[Brambilla {\em et~al.\/}, 2005]{Brambilla:2005}
    {A.~Brambilla, P.~Maffezzoni, \& G.~{Storti-Gajani}},
    ``Computation of period sensitivity functions for the simulation of phase noise in oscillators,''
    {\it {IEEE} Transactions on Circuits and Systems {I}: Regular Papers}, vol. 52, pp. 717--737, 2008.
%
\bibitem[Brambilla \& {Storti-Gajani}, 2008]{Brambilla:2008}
    {A.~Brambilla \& G.~{Storti-Gajani}},
    ``Computation of all the {Floquet} eigenfunctions in autonomous circuits,''
    {\it International Journal of Circuit Theory and Applications}, vol. 36, pp. 681--694, 2005.
%
\bibitem[Buonomo \& {Lo Schiavo}, 2003]{Buonomo:2003}
    {A.~Buonomo \& A.~{Lo Schiavo}},
    ``A constructive method for finding the periodic response of nonlinear circuits,''
    {\it {IEEE} Transactions on Circuits and Systems {I}: Fundamental Theory and Applications}, vol. 50, pp. 885--893, 2003.
%
\bibitem[Champneys \& Sandstede, 2007]{Champneys:2007}
    {{A.R.} Champneys \& B.~Sandstede},
    ``Numerical computation of coherent structures,''
    in {\it Numerical Continuation Methods for Dynamical Systems (B.~Krauskopf, {H.M.} Osinga and J.~Galan-Vioque, eds.)},
    Springer, pp. 331--358, 2007.
%
\bibitem[Cochelin \& Vergez, 2009]{Cochelin:2009}
    {B.~Cochelin \& C.~Vergez},
    ``A high order purely frequency-based harmonic balance formulation for continuation of periodic solutions,''
    {\it Journal of Sound and Vibration}, vol. 324, pp. 243--262, 2009.
%
\bibitem[Collado \& S\'uarez, 2005]{Collado:2005}
    {A.~Collado \& A.~S\'uarez},
    ``Application of bifurcation control to practical circuit design,''
    {\it IEEE Transactions on Microwave Theory and Techniques}, vol. 53, pp. 2777--2788, 2005.
%
\bibitem[{De Feo} {\em et~al.\/}, 2000]{DeFeoetal:2000}
    {O.~{De Feo}, {G.M.} Maggio, \& {M.P.} Kennedy},
    ``The {Colpitts} oscillator: {F}amilies of periodic solutions and their bifurcations,''
    {\it International Journal of Bifurcation and Chaos}, vol. 10, pp. 935--958, 2000.
%
\bibitem[{De Feo} \& Maggio, 2003]{DeFeoMaggio:2003}
    {O.~{De Feo} \& {G.M.} Maggio},
    ``Bifurcations in the {Colpitts} oscillator: {F}rom theory to practice,''
    {\it International Journal of Bifurcation and Chaos}, vol. 13, pp. 2917--2934, 2003.
%
\bibitem[Desroches {\em et~al.\/}, 2008]{Desroches:2008}
    {{M.F.} Desroches, B.~Krauskopf, \& {H.M.} Osinga},
    ``Mixed-mode oscillations and slow manifolds in the self-coupled FitzHugh Nagumo system,''
    {\it CHAOS}, vol. 18, p. 015107, 2008.
%
\bibitem[Dhooge {\em et~al.\/}, 2003]{Dhooge:2003}
    {A.~Dhooge, W.~Govaerts, \& {Y.A.} Kuznetsov},
    ``{MATCONT}: {A} {MATLAB} package for numerical bifurcation analysis of {ODE}s,''
    {\it ACM Trans. Math. Software}, vol. 29, pp. 141--164, 2003.
%
\bibitem[Doedel {\em et~al.\/}, 2006]{Doedel:2006}
    {{E.J.} Doedel, B.~Krauskopf, \& {H.M.} Osinga},
    ``Global bifurcations of the Lorenz manifold,''
    {\it Nonlinearity}, vol. 19, pp. 2947--2973, 2006.
%
\bibitem[Doedel \& Oldeman, 2009]{Auto07PMan}
    {{E.J.} Doedel \& {B.E.} Oldeman},
    ``{AUTO-07P}: {C}ontinuation and {B}ifurcation {S}oftware for {O}rdinary {D}ifferential {E}quations,''
    {C}oncordia {U}niversity, Montreal, Quebec, Canada, 2009.
%
\bibitem[Genesio {\em et~al.\/}, 1993]{Genesio:1993}
    {R.~Genesio, A.~Tesi, \& F.~Villoresi},
    ``A frequency approach for analyzing and controlling chaos in nonlinear circuits,''
    {\it {IEEE} Transactions on Circuits and Systems {I}: Fundamental Theory and Applications}, vol. 40, pp. 819--828, 1993.
%
\bibitem[Ghahramani {\em et~al.\/}, 2007]{Ghahramani:2007}
    {{M.M.} Ghahramani, S.~{Danseshgar Asl}, {M.P.} Kennedy, \& O.~{De Feo}},
    ``Optimizing the design of an injection-locked frequency divider by means of nonlinear analysis,''
    in {\it Proceedings of the European Conference on Circuit Theory and Design (ECCTD'07)}, pp. 571--574, August 2007.
%
\bibitem[Gourary {\em et~al.\/}, 2000]{Gourary:2000}
    {M.~Gourary, S.~Ulyanov, M.~Zharov, S.~Rusakov, {K.K.} Gullapalli, \& {B.J} Mulvaney},
    ``A robust and efficient oscillator analysis technique using harmonic balance,''
    {\it Computer Methods in Applied Mechanics and Engineering}, vol. 181, pp. 451--466, 2000.
%
\bibitem[Hu {\em et~al.\/}, 1989]{Hu:1989}
    {Y.~Hu, {J.J.} Obregon, \& {J.C.} Mollier},
    ``Nonlinear analysis of microwave {FET} oscillators using {V}olterra series,''
    {\it IEEE Transactions on Microwave Theory and Techniques}, vol. 37, pp. 1689--1693, 1989.
%
\bibitem[Huang \& Chu, 1994]{Huang:1994}
    {{C.C.} Huang \& {T.H.} Chu},
    ``Analysis of {MESFET} injection-locked oscillators in fundamental mode of operation,''
    {\it IEEE Transactions on Microwave Theory and Techniques}, vol. 42, pp. 1851--1857, 1994.
%
\bibitem[Kuznetsov, 2004]{Kuznetsov:2004}
    {Y.A. Kuznetsov},
    ``Elements of Applied Bifurcation Theory,''
    third ed., Springer, New York, 2004.
%
\bibitem[Lang {\em et~al.\/}, 2006]{Lang:2006}
    {{Z.Q.} Lang, {S.A.} Billings, {G.R.} Tomlinson, \& R.~Yue},
    ``Analytical description of the effects of system nonlinearities on output frequency responses: A case study,''
    {\it Journal of Sound and Vibration}, vol. 295, pp. 584--601, 2006.
%
\bibitem[Maggio {\em et~al.\/}, 1999]{Maggioetal:1999}
    {{G.M.} Maggio, O.~{De Feo}, \& {M.P.} Kennedy},
    ``Nonlinear analysis of the {Colpitts} oscillator and applications to design,''
    {\it {IEEE} Transactions on Circuits and Systems {I}: Fundamental Theory and Applications}, vol. 46, pp. 1118--1130, 1999.
%
\bibitem[Peng {\em et~al.\/}, 2008]{Peng:2008}
    {{Z.K.} Peng, {Z.Q.} Lang, {S.A.} Billings, \& {G.R.} Tomlinson},
    ``Comparisons between harmonic balance and nonlinear output frequency response function in nonlinear system analysis,''
    {\it Journal of Sound and Vibration}, vol. 311, pp. 56--73, 2008.
%
\bibitem[Rizzoli \& Neri, 2009]{Rizzoli:1988}
    {V.~Rizzoli \& A.~Neri},
    ``State of the art and present trends in nonlinear microwave {CAD} techniques,''
    {\it {IEEE} Transactions on Microwave Theory and Techniques}, vol. 36, pp. 343--364, 1998.
%
\bibitem[Schmidt \& Tondl, 1986]{Schmidt:1986}
    {G.~Schmidt \& A.~Tondl},
    ``Non-linear Vibrations,''
    Cambridge University Press, Cambridge, 1986.
%
\bibitem[Suarez {\em et~al.\/}, 2006]{Suarez:1998}
    {A.~S\'uarez, J.~Morales, \& R.~Qu\'er\'e},
    ``Synchronization analysis of autonomous microwave circuits using new global stability analysis tools,''
    {\it {IEEE} Transactions on Microwave Theory and Techniques}, vol. 46, pp. 494--504, 1998.
%
\end{thebibliography}
\end{document}